# Congruences between Selmer groups [*]

Li Guo

(liguo@newark.rutgers.edu)

Department of Mathematics and Computer Science

Rutgers University at Newark

Newark, NJ 07102

September 11, 1998

The study of congruences between arithmetically interesting numbers has a long history and plays important roles in several areas of number theory. Examples of such congruences include the Kummer congruences between Bernoulli numbers and congruences between coefficients of modular forms. Many of these congruences could be interpreted as congruences between special values of $L$-functions of arithmetic objects (motives). In recent years general conjectures [**B-K**] have been formulated relating these special values to the Selmer groups and other arithmetic invariants of the associated motives. In view of these conjectures, congruences between special values should give certain congruences between the corresponding Selmer groups.

In this paper, we take a different point of view by studying congruences between Selmer groups and deducing consequences of such congruences to the congruences of special values. Roughly speaking, given two Galois representations that are congruent on a finite level, i.e., that become isomorphic representations modulo a prime power, we consider the relation between the corresponding Selmer groups. Precise definitions will be given in later sections. We will focus on the Selmer groups defined by Bloch and Kato [**B-K**] and will provide congruences when the congruent Galois representations are from cyclotomic characters, Hecke characters from CM elliptic curves, and from adjoints of modular forms. We also display consequences of such congruences to the congruences of special values of $L$-functions, in particular, to special values of the Riemann zeta function. It is interesting to observe that the congruences of special values obtained this way are different from the classical congruences, such as the Kummer congruences.

Methods used in this paper are mostly from Iwasawa theory, including the classical theory originated from Iwasawa [**Iw**] and the "horizontal" theory recently developed

---

[*]This research is supported in part by NSF grant #DMS 97-96122. MSC Number: Primary 11R23. Secondary 14G10, 14G25.

by Wiles [**Wi**], even though only the algebraic part of the theory is applied.

There are several versions of Selmer groups, motivated by the Selmer groups of elliptic curves. We choose to study the Selmer group of Bloch and Kato since it is the one that applies to the most general context and relates well with spcial values. However, in order to get a congruence, we have to make use of other versions of Selmer groups, defined by Greenberg [**Gr**] and by Wiles [**Wi**]. This suggests that the Selmer groups of Greenberg and Wiles might be better candidates for studying congruence problems. For some other results on congruences of Selmer groups, see [**Gr**].

The layout of this paper is as follows. In §1, congruences between Galois representations and between abelian groups are defined. We introduce several types of congruences and study the relation between them (Proposition 1.1). For later applications, we also give general criteria for congruences of Selmer groups (Proposition 1.2). In §2, we prove congruences (Theorem 1) between Selmer groups from cyclotomic representations. By the theorem of Bloch and Kato that relates these Selmer groups to special values of the Riemann zeta function, we are able to obtain congruences between these special value, and congruences between the Bernoulli numbers (Theorem 2). As indicated earlier, these congruences are different from the classical Kummer congruences. We then prove, in §3, congruences when the Selmer groups are from Hecke characters over an imaginary quadratic field $K$. This could be regarded as an algebraic analogue of the congruences between special values of Hecke $L$-functions that were used to construct $p$-adic $L$-functions. The result on primes that split in $K$ is quite complete (Theorem 3). However, for primes non-split in $K$, we have to impose an additional condition, called r-congruence, on the Galois representations in order to obtain congruences between the Selmer groups (Theorem 5). Congruences in §2 and §3 should still hold when the characters are twisted by a Dirichlet character. In the last section, we provide congruences between adjoints of two dimensional Galois representations (Theorem 6). As a special case, we obtain congruences between adjoint representations from modular forms (Corollary 4.1). Unlike the cases of cyclotomic characters and Hecke characters, congruences between special values from the adjoint representations are not know before. However, with the help of the Bloch-Kato conjecture, congruences from Corollary 4.1 give congruences between the corresponding special values.

The author would like to thank F. Diamond, R. Greenberg, K. Rubin and J. Sturm for helpful discussions.

# 1 Definitions and preliminary results

In this section we recall the definition of Selmer groups of Bloch and Kato, define the congruence between two Selmer groups and give preliminary results on congruences between Selmer groups. For the rest of this paper, we denote $F$ and $K$ for number fields, denote $\mathcal{O}_F$ and $\mathcal{O}_K$ for rings of integers, denote $\mathfrak{p}$ (resp. $\lambda$) for a prime of $F$ (resp. of $K$), and denote $\mathcal{O}_\mathfrak{p}$ and $\mathcal{O}_\lambda$ for the rings of integers of the completions $F_\mathfrak{p}$ and $K_\lambda$.

## 1.1 Selmer groups

Given a finite prime $\lambda$ of $K$, by a $\lambda$-adic represetation of $G_F$, or a $\mathcal{O}_\lambda[[G_F]]$-module, we mean a topological $\mathcal{O}_\lambda$-module $A$ with a continuous action of $G_F$. The representations we consider are of the following five types.

1. a finite dimensional $F_\lambda$ vector space $A$ with the product topology and a continuous group homomorphism $\rho: G_F \to \mathrm{Aut}_{\mathrm{cont}}(A)$;

2. a finitely generated, free $\mathcal{O}_\lambda$-module $A$ with the compact topology and a continuous group homomorphism $\rho: G_F \to \mathrm{Aut}_{\mathrm{cont}}(A)$;

3. a cofinitely generated, cofree $\mathcal{O}_\lambda$-module $A$ (i.e., $A \cong (K_\lambda/\mathcal{O}_\lambda)^d$) with the discrete topology and a continuous group homomorphism $\rho: G_F \to \mathrm{Aut}_{\mathrm{cont}}(A)$. Such a representation can be identified with the Pontrjagin dual of a representation of type (2) with the dual action of $G_F$. The integer $d$ is called the corank of $A$, denoted by $\mathrm{corank}_{\mathcal{O}_\lambda} A$;

4. a quotient $A/\lambda^n A$, $n \geq 1$ for the $A$ of type (2) with the quotient action of $G_F$;

5. a submodule $A[\lambda^n] \stackrel{\mathrm{def}}{=} \{x \in A \mid \lambda^n x = 0\}$ for the $A$ of type (3) with the restricted action of $G_F$.

We further assume that the action of $G_F$ on $A$ is unramified outside of a finite set of primes $S$ of $F$ that contains $\infty$ and primes of $F$ above $\ell$. Then the action of $G_F$ is naturally induced from the action of $G_S$, the Galois group of the maximal extension of $F$ unramified outside of $S$.

Let $A$ be a representation of $G_F$ of one of the five types listed above. A local condition $L$ for $A$ consists of the selection of a subgroup

$$H^1_L(F_\mathfrak{p}, A) \subseteq H^1(F_\mathfrak{p}, A)$$

for each finite prime $\mathfrak{p}$ of $F$. The Selmer group of $A$ with respect to the local condition $L$ is defined by

$$H^1_L(F, A) = \ker\left\{H^1(F, A) \to \bigoplus_\mathfrak{p} \frac{H^1(F_\mathfrak{p}, A)}{H^1_L(F_\mathfrak{p}, A)}\right\}.$$

The Selmer group that we will study in this paper is defined by the local condition of Bloch and Kato [**B-K**], for a representation $A$ of type (3). Given such an $A$, define

$$T = \varprojlim{}_n A[\lambda^n] \text{ and } V = T \otimes_{\mathcal{O}_\lambda} K_\lambda$$

with the induced representations of $G_F$. Then $T$ is of type (2) and $V$ is of type (1) and $A \cong V/T$ as representations of $G_F$. For each finite prime $\mathfrak{p}$ of $F$, let $p$ be the prime of $\mathbb{Q}$ under $\mathfrak{p}$ and let $B_{\text{crys},\mathfrak{p}}$ be the ring defined by Fontaine. Define

$$H^1_f(F_\mathfrak{p}, V) = \begin{cases} H^1(G_\mathfrak{p}/I_\mathfrak{p}, V^{I_\mathfrak{p}}), & \text{if } \lambda \nmid p \\ \ker\{H^1(F_\mathfrak{p}, V) \to H^1(F_\mathfrak{p}, B_{\text{crys},\mathfrak{p}} \otimes V)\}, & \text{if } \lambda \mid p \end{cases}$$

where $G_\mathfrak{p}$ and $I_\mathfrak{p}$ are the decomposition group and initial group of $\mathfrak{p}$ with respect to a fixed embedding $\bar{F} \to \bar{F}_\mathfrak{p}$. Let $\text{pr} : V \to A$ be the natural projection. Define

$$H^1_f(F_\mathfrak{p}, A) = \text{pr}(H^1_f(F_\mathfrak{p}, V)),$$

and define the Selmer group of Bloch and Kato by

$$H^1_f(F, A) = \ker\left\{H^1(F, A) \to \bigoplus_\mathfrak{p} \frac{H^1(F_\mathfrak{p}, A)}{H^1_f(F_\mathfrak{p}, A)}\right\}.$$

It is clear that $H^1_f(F, A)$ is a torsion abelian group.

## 1.2 Congruences

We now define the congruences between two torsion abelian groups. This definition applies in particular to the Selmer groups. Let $Y$ be a finite set. Denote $|Y|$ for the cardinality of $Y$.

**Definition 1.1** *Let $M_1$ and $M_2$ be two torsion $\mathcal{O}_K$-modules.*

1. *$M_1$ and $M_2$ are called **algebraically congruent** modulo $\lambda^n$, denoted $M_1 \equiv_{alg} M_2 \bmod \lambda^n$, if $M_1[\lambda^n] \cong M_2[\lambda^n]$;*

2. *$M_1$ and $M_2$ are called **numerically congruent** modulo $\lambda^n$, denoted $M_1 \equiv_{num} M_2 \bmod \lambda^n$, if $|M_1[\lambda^n]|$ and $|M_2[\lambda^n]|$ are finite and $|M_1[\lambda^n]|=|M_2[\lambda^n]|$;*

3. *$M_1$ and $M_2$ are **cardinally congruent** modulo $\lambda^n$, denoted $M_1 \equiv_{car} M_2 \bmod \lambda^n$, if $M_1$ and $M_2$ are finite and $|M_1|\equiv|M_2| \bmod \lambda^n$.*

We first display the following simple relations between the three types of congruences.

**Proposition 1.1** *Let $\lambda$ be a finite prime of $K$, and let $M_1$ and $M_2$ be two finite $\mathcal{O}_\lambda$-modules.*

*$M_1$ and $M_2$ are algebraically congruent modulo $\lambda^n$ if and only if $M_1$ and $M_2$ are numerically congruent modulo $\lambda^r$, for all $1 \leq r \leq n$;*

2. If $M_1$ and $M_2$ are numerically congruent modulo $\lambda^n$, then $M_1$ and $M_2$ are cardinally congruent modulo $\lambda^n$.

**Proof:** 1. By definition, $M_1 \equiv_{alg} M_2 \mod \lambda^n$ implies that $M_1[\lambda^r] \cong M_2[\lambda^r]$ for $r = 1, \cdots, n$. This shows that $M_1 \equiv_{num} M_2 \mod \lambda^i$ for $r = 1, \cdots, n$. To prove the implication in the other direction, we apply induction on $n$. Since $\mathcal{O}_\lambda$ is a PID, we could write, for $i = 1, 2$,
$$M_i = \oplus_{j=1}^{k_i} \mathcal{O}_\lambda / \lambda^{c_{i,j}}.$$

For $i = 1, 2$ and for each integer $r \geq 1$, denote $s_{i,r} = |\{a_{i,j}, 1 \leq j \leq k | a_{i,j} \geq r\}|$ and $t_{i,r} = |\{a_{i,j}, 1 \leq j \leq k | a_{i,j} = r\}|$. Then $t_{i,r} = s_{i,r} - s_{i,r+1}$. Also
$$M_i[\lambda^r] = \oplus_{u=1}^{r-1} (\mathcal{O}_\lambda / \lambda^u)^{t_{i,u}} \bigoplus (\mathcal{O}_\lambda / \lambda^r)^{s_{i,r}}.$$

Therefore $M_1[\lambda^r] \cong M_2[\lambda^r]$ if and only if $t_{1,u} = t_{2,u}$ for $1 \leq u \leq r - 1$ and $s_{1,r} = s_{2,r}$. On the other hand,
$$|M_i[\lambda^r]| = |\mathcal{O}_\lambda / \lambda|^{\sum_{u=1}^{r-1} u t_{i,u} + r s_{i,r}}.$$

When $n = 1$, $M_1$ and $M_2$ being congruent module $\lambda^n$ means $|M_1[\lambda]| = |M_2[\lambda]|$. So $s_{1,1} = s_{2,1}$. Therefore $M_1[\lambda] \cong M_2[\lambda]$. Now inductively assume that the converse is true for $n$, and assume that $|M_1[\lambda^r]| = |M_1[\lambda^r]|$ for $r = 1, \cdots n + 1$. By induction, we have $M_1[\lambda^n] \cong M_2[\lambda^n]$. Then

$$t_{1,u} = t_{2,u}, \quad 1 \leq u \leq n-1, \text{ and } s_{1,n} = s_{2,n}. \tag{1}$$

From $|M_1[\lambda^{n+1}]| = |M_2[\lambda^{n+1}]|$ we have

$$\sum_{u=1}^{n} u t_{1,u} + (n+1) s_{1,n+1} = \sum_{u=1}^{n} u t_{2,u} + (n+1) s_{2,n+1}. \tag{2}$$

Solving Eq (1) and Eq (2), observing that $s_{i,n} = t_{i,n} + s_{i,n+1}$, we obtain

$$s_{1,n+1} = s_{2,n+1}, \ t_{1,n} = t_{2,n}.$$

Together with Eq (1), this implies $M_1[\lambda^{n+1}] \cong M_2[\lambda^{n+1}]$. This completes the induction.

2. If $M_1 \equiv_{num} M_2 \mod \lambda^n$, then by definition, $|M_1[\lambda^n]| = |M_2[\lambda^n]|$. Consider the exact sequence of $\mathcal{O}_\lambda$-modules

$$0 \to M_i[\lambda^n] \to M_i \xrightarrow{\lambda^n} \lambda^n M_i \to 0, \tag{3}$$

where we identify $\lambda$ with a generator of the principle maximal ideal $\lambda$ of $\mathcal{O}_\lambda$. First consider the case when $\lambda^n M_1 \neq 0$. Let $x \in M_1$ with $\lambda^n x \neq 0$. Then $(\mathcal{O}_\lambda x)[\lambda^n] \cong \mathcal{O}_\lambda / \lambda^n$. So $|M_1[\lambda^n]| \equiv 0 \pmod{\lambda^n}$. Therefore $\lambda^n$ divides $|M_2[\lambda^n]|$. Thus we have

$$|M_1| \equiv |M_2| (\equiv 0) \mod \lambda^n.$$

The case when $\lambda^n M_2 \neq 0$ can be similarly proved. Now consider the case when $\lambda^n M_1 = 0$ and $\lambda^n M_2 = 0$. Then from the exact sequence (3), we must have $M_i = M_i[\lambda^n]$, $i = 1, 2$. Therefore, $|M_1[\lambda^n]| = |M_2[\lambda^n]|$ means that $|M_1| = |M_2|$. This of course implies $|M_1| \cong |M_2| \mod \lambda^n$. ∎

**Remark:** Note that $M_1 \equiv_{car} M_2 \mod \lambda^n$ does not imply $M_1 \equiv_{num} M_2 \mod \lambda^n$. For example, let $M_1 = \mathbb{Z}/\ell^5\mathbb{Z}$, $M_2 = \mathbb{Z}/\ell^4\mathbb{Z} \oplus \mathbb{Z}/\ell\mathbb{Z}$. Then $M_1 \equiv_{car} M_2 \mod \ell$, but $M_1 \not\equiv_{num} M_2 \mod \ell$.

Let $A_1$ and $A_2$ be two $\mathcal{O}_\lambda[[G_F]]$-modules that are cofree $\mathcal{O}_\lambda$-modules with finite coranks over $\mathcal{O}_\lambda$. Suppose $A_1[\lambda^n]$ and $A_2[\lambda^n]$ are isomorphic $(\mathcal{O}_\lambda/\lambda^n)[[G_F]]$-modules. We would like to study the conditions under which the corresponding Selmer groups $H^1_f(F, A_1)$ and $H^1_f(F, A_2)$ are (algebraically) congruent modulo $\lambda^n$. We will provide some criteria in Proposition 1.2. Applications will be given in later sections.

Let $A$ be either $A_1$ or $A_2$. We now define a Selmer group of $A$ on "level $r$" for each integer $r \geq 1$. For each finite prime $\mathfrak{p}$ of $F$, from the exact sequence

$$0 \to A[\lambda^r] \xrightarrow{j_r} A \xrightarrow{\lambda^r} A \to 0$$

we get the exact sequence

$$0 \to \frac{H^0(F_\mathfrak{p}, A)}{\lambda^r H^0(F_\mathfrak{p}, A)} \to H^1(F_\mathfrak{p}, A[\lambda^r]) \xrightarrow{j_r} H^1(F_\mathfrak{p}, A)[\lambda^r] \to 0.$$

Define the subgroup of $H^1(F_\mathfrak{p}, A[\lambda^n])$

$$H^1_f(F_\mathfrak{p}, A, \lambda^r) = j_r^{-1}(H^1_f(F_\mathfrak{p}, A)[\lambda^r]),$$

and define the subgroup of $H^1(F, A[\lambda^n])$

$$H^1_f(F, A, \lambda^r) = \ker\left\{H^1(F, A[\lambda^r]) \to \bigoplus_\mathfrak{p} \frac{H^1(F_\mathfrak{p}, A[\lambda^r])}{H^1_f(F_\mathfrak{p}, A, \lambda^r)}\right\}.$$

Note that $H^1_f(F, A, \lambda^r)$ depends on $A$, not just on $A[\lambda^r]$. In particular, even if $A_1[\lambda^r] \cong A_2[\lambda^r]$ as $G_F$-modules, it is not necessarily true that

$$H^1_f(F, A_1, \lambda^n) \cong H^1_f(F, A_2, \lambda^n).$$

**Proposition 1.2** *Let $A_1$ and $A_2$ be two $\mathcal{O}_\lambda[[G_F]]$-modules that are cofree over $\mathcal{O}_\lambda$ with finite corank. Assume that $A_1[\lambda^n]$ and $A_2[\lambda^n]$ are isomorphic as $(\mathcal{O}_\lambda/\lambda^n)[[G_F]]$-modules.*

  1. *If there is a $\mathcal{O}_\lambda[[G_F]]$-module $X$ such that $H^1_f(F, A_i) \cong \mathrm{Hom}_{G_F}(X, A_i)$ for $i = 1, 2$, then $H^1_f(F, A_1) \equiv_{alg} H^1_f(F, A_2) \pmod{\lambda^n}$.*

2. Assume that $\operatorname{corank}_{\mathcal{O}_\lambda} H^0(F, A_1) = \operatorname{corank}_{\mathcal{O}_\lambda} H^0(F, A_2)$. Then $H^1_f(F, A_1) \equiv_{alg} H^1_f(F, A_2) \pmod{\lambda^n}$ if and only if $|H^1_f(F, A_1, \lambda^r)| = |H^1_f(F, A_2, \lambda^r)|$ for $1 \le r \le n$.

3. Assume that $\operatorname{corank}_{\mathcal{O}_\lambda} H^0(F, A_1) = \operatorname{corank}_{\mathcal{O}_\lambda} H^0(F, A_2)$. Fix an $\mathcal{O}_\lambda[[G_F]]$-module isomorphism $\sigma : A_1[\lambda^n] \to A_2[\lambda^n]$ and for each finite prime $\mathfrak{p}$ of $F$ and each $1 \le r \le n$, let $\sigma_{\mathfrak{p},r} : H^1(F_\mathfrak{p}, A_1[\lambda^r]) \to H^1(F_\mathfrak{p}, A_2[\lambda^r])$ be the induced isomorphism. If $\sigma_{\mathfrak{p},r}(H^1_f(F_\mathfrak{p}, A_1, \lambda^r)) = H^1_f(F_\mathfrak{p}, A_2, \lambda^r)$ for all finite prime $\mathfrak{p}$ and $1 \le r \le n$, then $H^1_f(F, A_1) \equiv_{alg} H^1_f(F, A_2) \pmod{\lambda^n}$.

**Proof:** 1. For $i = 1, 2$, from the exact sequence

$$0 \to A_i[\lambda^n] \to A_i \xrightarrow{\lambda^n} A_i \to 0$$

we get the exact sequence

$$0 \to \operatorname{Hom}_{G_F}(X, A_i[\lambda^n]) \to \operatorname{Hom}_{G_F}(X, A_i) \xrightarrow{\lambda^n} \operatorname{Hom}_{G_F}(X, A_i).$$

Thus we have the isomorphism of abelian groups

$$\operatorname{Hom}_{G_F}(X, A_i[\lambda^n]) \cong \operatorname{Hom}_{G_F}(X, A_i)[\lambda^n].$$

Therefore from the assumptions we have

$$\begin{aligned} H^1_f(F, A_1)[\lambda^n] &\cong \operatorname{Hom}_{G_F}(X, A_1)[\lambda^n] \\ &\cong \operatorname{Hom}_{G_F}(X, A_1[\lambda^n]) \\ &\cong \operatorname{Hom}_{G_F}(X, A_2[\lambda^n]) \\ &\cong \operatorname{Hom}_{G_F}(X, A_2)[\lambda^n] \\ &\cong H^1_f(F, A_2)[\lambda^n]. \end{aligned}$$

2. Let $A$ be $A_1$ or $A_2$. By definition, we have the commutative diagram with exact rows

$$\begin{array}{ccccccccc} 0 & \to & \frac{H^0(F,A)}{\lambda^r H^0(F,A)} & \to & H^1(F, A[\lambda^r]) & \to & H^1(F, A)[\lambda^r] & \to & 0 \\ & & \downarrow & & \downarrow & & \downarrow & & \\ 0 & \to & 0 & \to & \bigoplus_\mathfrak{p} \frac{H^1(F_\mathfrak{p}, A[\lambda^r])}{H^1_f(F_\mathfrak{p}, A, \lambda^r)} & \to & \bigoplus_\mathfrak{p} \frac{H^1(F_\mathfrak{p}, A)}{H^1_f(F_\mathfrak{p}, A)} & & \end{array}$$

Applying the snake lemma we obtain the exact sequence

$$0 \to \frac{H^0(F, A)}{\lambda^r H^0(F, A)} \to H^1_f(F, A, \lambda^r) \xrightarrow{j_r} H^1_f(F, A)[\lambda^r] \to 0. \tag{4}$$

Let $H_{\text{div}}$ be the maximal $\lambda$-divisible subgroup of $H^0(F, A)$ and let $H_{\text{cot}}$ be the cotorsion quotient $H^0(F, A)/H_{\text{div}}$. Then we have the exact sequence

$$0 \to H_{\mathrm{div}} \to H^0(F, A) \to H_{\mathrm{cot}} \to 0 \tag{5}$$

of $\mathcal{O}_\lambda$-modules. Let $d$ be the $\mathcal{O}_\lambda$-corank of $H^0(F, A)$, then $H_{\mathrm{div}} \cong (K_\lambda/\mathcal{O}_\lambda)^d$. Since $H_{\mathrm{cot}}$ is finite, from the exact sequence

$$0 \to H_{\mathrm{cot}}[\lambda^r] \to H_{\mathrm{cot}} \xrightarrow{\lambda^r} H_{\mathrm{cot}} \to \frac{H_{\mathrm{cot}}}{\lambda^r H_{\mathrm{cot}}} \to 0$$

we get

$$|H_{\mathrm{cot}}[\lambda^r]| = |\frac{H_{\mathrm{cot}}}{\lambda^n H_{\mathrm{cot}}}| . \tag{6}$$

Multiplying the sequence (5) by $\lambda^r$ and taking the snake lemma exact sequence from the resulting commutative diagram, we obtain the exact sequence

$$0 \to (\mathcal{O}_\lambda/\lambda^r)^d \to H^0(F, A)[\lambda^r] \to H_{\mathrm{cot}}[\lambda^r]$$
$$\to 0 \to \frac{H^0(F, A)}{\lambda^r H^0(F, A)} \to \frac{H_{\mathrm{cot}}}{\lambda^r H_{\mathrm{cot}}} \to 0.$$

This, together with equation (6) and the equation $H^0(F, A[\lambda^r]) = H^0(F, A)[\lambda^r]$, enables us to get

$$\begin{aligned}
|\frac{H^0(F, A)}{\lambda^r H^0(F, A)}| &= |\frac{H_{\mathrm{cot}}}{\lambda^r H_{\mathrm{cot}}}| \\
&= |H_{\mathrm{cot}}[\lambda^r]| \\
&= |H^0(F, A[\lambda^r])| \, |\mathcal{O}_\lambda/\lambda^r|^{-d} .
\end{aligned}$$

Combining this with the sequence (4), we obtain

$$|H^1_f(F, A)[\lambda^r]| = |H^1_f(F, A, \lambda^r)| \, |H^0(F, A[\lambda^r])|^{-1} \, |\mathcal{O}_\lambda/\lambda^r|^d .$$

Now the assumption that $\mathrm{corank}_{\mathcal{O}_\lambda} H^0(F, A_1) = \mathrm{corank}_{\mathcal{O}_\lambda} H^0(F, A_2)$ implies $d_1 = d_2$, and the assumption that $A_1[\lambda^r] \cong A_2[\lambda^r]$ implies $|H^0(F, A_1[\lambda^r])| = |H^0(F, A_2[\lambda^r])|$ . Thus the statement that $|H^1(F, A_1)[\lambda^r]| = |H^1(F, A_1)[\lambda^r]|$ for $1 \leq r \leq n$ is equivalent to the statement that $|H^1(F, A_1, \lambda^r)| = |H^1(F, A_1, \lambda^r)|$ for $1 \leq r \leq n$. By Proposition 1.1, the first statement is true if and only $H^1_f(F, A_1)[\lambda^n] \cong H^1_f(F, A_2)[\lambda^n]$.

3. Fix an $r$ with $1 \leq r \leq n$. Let $\sigma_r : H^1(F, A_1[\lambda^r]) \to H^1(F, A_2[\lambda^r])$ be the isomorphism induced by $\sigma$. For each finite prime $\mathfrak{p}$ of $F$, we have the commutative diagram

$$\begin{array}{ccc}
H^1(F, A_1[\lambda^r]) & \xrightarrow{\sigma_r} & H^1(F, A_2[\lambda^r]) \\
\downarrow & & \downarrow \\
H^1(F_\mathfrak{p}, A_1[\lambda^r]) & \xrightarrow{\sigma_{\mathfrak{p},r}} & H^1(F_\mathfrak{p}, A_2[\lambda^r])
\end{array}$$

where the vertical maps are restriction maps. Then by the definition of $H^1_f(F, A_i, \lambda^r)$ with $i = 1, 2$ we have $\sigma_r(H^1_f(F, A_1, \lambda^n)) = H^1_f(F, A_2, \lambda^r)$. In particular, we have $|H^1_f(F, A_1, \lambda^n)| = |H^1_f(F, A_2, \lambda^r)|$. It now follows from part 2 of the proposition that $H^1_f(F, A_1) \equiv_{alg} H^1_f(F, A_2) \pmod{\lambda^n}$. ∎

## 2 Cyclotomic characters and Kummer congruences

In this section we first study the congruences of Selmer groups from the cyclotomic characters. We then compare these congruences with the Kummer congruences.

### 2.1 Algebraic congruences

Let $\ell$ be a fixed prime. Let $\chi = \chi_\ell : G_\mathbb{Q} \to \mathbb{Z}_\ell^\times$ be the $\ell$-adic cyclotomic character. Thus $\chi_\ell$ defines a one dimensional Galois representation which is denoted $\mathbb{Z}_\ell(1)$. For any integer $n$, the character $\chi_\ell^n$ defines a one dimensional representation $\mathbb{Z}_\ell(n) = \mathbb{Z}_\ell(1)^{\otimes n}$ and the induced representation

$$(\mathbb{Q}_\ell/\mathbb{Z}_\ell)(n) \stackrel{\text{def}}{=} \mathbb{Z}_\ell(n) \otimes (\mathbb{Q}_\ell/\mathbb{Z}_\ell).$$

**Theorem 1** *Let $\ell$ be an odd prime and let $n$ be a positive integer. If $k$ and $k'$ are two integers both greater then 1 or both less then 0 such that $\chi^k \equiv \chi^{k'} \bmod \ell^n$, then*

$$H^1_f(\mathbb{Q}, (\mathbb{Q}_\ell/\mathbb{Z}_\ell)(\chi^k)) \equiv_{alg} H^1_f(\mathbb{Q}, (\mathbb{Q}_\ell/\mathbb{Z}_\ell)(\chi^{k'})) \bmod \ell^n.$$

**Proof:** The following proof is suggested by Rubin. One could also use Iwasawa theory, as in Theorem 3, and apply generalized Cassels-Tate pairing [**Fl**, **Gu1**].

To simplify notations, denote $A_1 = (\mathbb{Q}_\ell/\mathbb{Z}_\ell)(k)$ and $A_2 = (\mathbb{Q}_\ell/\mathbb{Z}_\ell)(k')$. By assumption, $A_1[\ell^n]$ is isomorphic to $A_2[\ell^n]$ as $\mathbb{Z}_\ell[[G_\mathbb{Q}]]$-modules. Fix an $\mathbb{Z}_\ell[[G_\mathbb{Q}]]$-module isomorphism $\sigma : A_1[\ell^n] \to A_2[\ell^n]$. For each finite prime $p$ of $\mathbb{Q}$ and each $1 \leq r \leq n$, let $\sigma_{p,r} : H^1(\mathbb{Q}_p, A_1[\ell^r]) \to H^1(\mathbb{Q}_p, A_2[\ell^r])$ be the induced isomorphism. We will prove that $\sigma_{p,r}(H^1_f(\mathbb{Q}_p, A_1, \ell^r)) = H^1_f(\mathbb{Q}_p, A_2, \ell^r)$ for all finite prime $\mathfrak{p}$ and $1 \leq r \leq n$. Then the proof follows from Proposition 1.2.

We first consider the case when $p \neq \ell$. Then both $A_1$ and $A_2$ are unramified at $p$. So for $i = 1, 2$, $H^1_f(\mathbb{Q}_p, A_i) = H^1(g_p, A_i)$. We also have the short exact sequence

$$0 \to A_i[\ell^r] \to A_i \xrightarrow{\ell^r} A_i \to 0$$

of $g_p$-modules. This, together with the inflation map, gives the following commutative diagram with exact rows

$$\begin{array}{ccccccccc}
0 & \to & H^0(g_p, A_i)/\ell^r H^0(g_p, A_i) & \to & H^1(g_p, A_i[\ell^r]) & \to & H^1(g_p, A_i)[\ell^r] & \to & 0 \\
& & \| & & \downarrow & & \downarrow & & \\
0 & \to & H^0(\mathbb{Q}_p, A_i)/\ell^r H^0(\mathbb{Q}_p, A_i) & \to & H^1(\mathbb{Q}_p, A_i[\ell^r]) & \to & H^1(\mathbb{Q}_p, A_i)[\ell^r] & \to & 0
\end{array}$$

in which the vertical maps are injective. This diagram shows that $H^1_f(\mathbb{Q}_p, A_i[\ell^r]) = H^1(g_p, A_i[\ell^r])$. It follows that $\sigma_{p,r}(H^1_f(\mathbb{Q}_p, A_1[\ell^n])) = (H^1_f(\mathbb{Q}_p, A_2[\ell^n]))$.

We now consider the case when $p = \ell$. First assume that $k, k' < 0$. By Example 3.9 from [**B-K**], $H^1_f(\mathbb{Q}_\ell, \mathbb{Q}_\ell(k)) = H^1_f(\mathbb{Q}_\ell, \mathbb{Q}_\ell(k')) = 0$. Then $H^1_f(\mathbb{Q}_\ell, A_i) = 0$ for $i = 1, 2$. Therefore $H^1_f(\mathbb{Q}_\ell, A_i, \ell^r) = \ker\{H^1(\mathbb{Q}_\ell, A_i[\ell^r]) \to H^1(\mathbb{Q}_\ell, A_i)\}$. Denote $L = \mathbb{Q}_\ell(\mu_\infty)$. Then the action of $\mathrm{Gal}(\bar{\mathbb{Q}}_\ell/\mathbb{Q}_\ell)$ on $A_i$ is induced from the action of $\mathrm{Gal}(L/\mathbb{Q}_\ell)$ and the inflation-restriction sequence gives us the following commutative diagram with exact rows

$$
\begin{array}{ccccccc}
0 & \to & H^1(\mathrm{Gal}(L/\mathbb{Q}_\ell), A_i[\ell^r]) & \to & H^1(\mathbb{Q}_\ell, A_i[\ell^r]) & \to & H^1(L, A_i[\ell^r]) \\
& & \downarrow & & \downarrow & & \downarrow \\
0 & \to & H^1(\mathrm{Gal}(L/\mathbb{Q}_\ell), A_i) & \to & H^1(\mathbb{Q}_\ell, A_i) & \to & H^1(L, A_i)
\end{array}
$$

Here the right vertical map is injective since $H^1(L, A_i[\ell^r]) = \mathrm{Hom}(G_L, A_i[\ell^r])$ and $H^1(L, A_i) = \mathrm{Hom}(G_L, A_i)$. Thus $H^1_f(\mathbb{Q}, A_i, \ell^r)$, defined as the kernel of the middle vertical map, in contained in the image of $H^1(\mathrm{Gal}(L/\mathbb{Q}_\ell), A_i[\ell^r])$ in $H^1(\mathbb{Q}_\ell, A_i[\ell^r])$. Since $\mathrm{Gal}(L/\mathbb{Q}_\ell)$ is topologically cyclic, we have

$$H^1(\mathrm{Gal}(L/\mathbb{Q}_\ell), A_i[\ell^r]) \cong A_i[\ell^r]/(\gamma - 1)A_i[\ell^r]$$

where $\gamma$ is a topological generator of $\mathrm{Gal}(L/\mathbb{Q}_\ell)$. This quotient is a cyclic abelian group. Therefore $H^1_f(\mathbb{Q}_\ell, A_i, \ell^r)$ is the unique subgroup of $H^1(\mathrm{Gal}(L/\mathbb{Q}_\ell), A_i[\ell^r])$ with order $\mid H^1_f(\mathbb{Q}_\ell, A_i, \ell^r) \mid$. From the exact sequence

$$0 \to H^0(\mathbb{Q}_\ell, A_i[\ell^r]) \to H^0(\mathbb{Q}_\ell, A_i) \xrightarrow{p^r} H^0(\mathbb{Q}_\ell, A_i) \to H^1_f(\mathbb{Q}_\ell, A_i, \ell^r) \to 0$$

of finite groups, we have $\mid H^1_f(\mathbb{Q}_\ell, A_i, \ell^r) \mid = \mid H^0(\mathbb{Q}_\ell, A_i[\ell^r]) \mid$. Now the isomorphism $\sigma : A_1[\ell^r] \to A_2[\ell^r]$ also induces an isomorphism from $H^1(\mathrm{Gal}(L/\mathbb{Q}_\ell), A_1[\ell^r])$ onto $H^1(\mathrm{Gal}(L/\mathbb{Q}_\ell), A_2[\ell^r])$ such that the resulting diagram

$$
\begin{array}{ccc}
H^1(\mathrm{Gal}(L/\mathbb{Q}_\ell), A_1[\ell^r]) & \hookrightarrow & H^1(\mathbb{Q}_\ell, A_1[\ell^r]) \\
\downarrow & & \downarrow \\
H^1(\mathrm{Gal}(L/\mathbb{Q}_\ell), A_2[\ell^r]) & \hookrightarrow & H^1(\mathbb{Q}_\ell, A_2[\ell^r])
\end{array}
$$

commutes. Therefore, the isomorphism $\sigma_{\ell,r} : H^1(\mathbb{Q}_\ell, A_1[\ell^r]) \cong H^1(\mathbb{Q}_\ell, A_2[\ell^r])$ sends $H^1_f(\mathbb{Q}_\ell, A_1, \ell^r)$, the unique subgroup of $H^1(\mathrm{Gal}(L/\mathbb{Q}_\ell), A_1[\ell^r]) \subseteq H^1(\mathbb{Q}_\ell, A_1[\ell^r])$ of order $H^0(\mathbb{Q}_\ell, A_1[\ell^r])$, onto $H^1_f(\mathbb{Q}_\ell, A_2, \ell^r)$, the unique subgroup of $H^1(\mathrm{Gal}(L/\mathbb{Q}_\ell), A_2[\ell^r]) \subseteq H^1(\mathbb{Q}_\ell, A_2[\ell^r])$ of order $\mid H^0(\mathbb{Q}_\ell, A_2[\ell^r]) \mid = \mid H^0(\mathbb{Q}_\ell, A_1[\ell^r]) \mid$. Now the theorem follows from Proposition 1.2, under the assumption $k, k' < 0$.

Next assume that $k, k' > 1$ with $\chi^k \equiv \chi^{k'} \pmod{\ell^n}$. Then $1 - k, 1 - k' < 0$ and $\chi^{1-k} \equiv \chi^{1-k'} \pmod{\ell^n}$. Denote $B_1$ and $B_2$ for the representations of $\mathbb{Z}_\ell$-corank one associated to $\mathbb{Q}_\ell(1-k)$ and $\mathbb{Q}_\ell(1-k')$, i.e., $B_1 = (\mathbb{Q}_\ell/\mathbb{Z}_\ell)(1-k)$, $B_2 = (\mathbb{Q}_\ell/\mathbb{Z}_\ell)(1-k')$. Let $\tau : B_1[\ell^n] \to B_2[\ell^n]$ be the isomorphism induced from $\sigma : A_1[\ell^n] \to A_2[\ell^n]$, and for each $r$ with $1 \leq r \leq n$, let $\tau_{\ell,r} : H^1(\mathbb{Q}_\ell, B_1[\ell^r]) \to H^1(\mathbb{Q}_\ell, B_2[\ell^r])$ be the isomorphism induced by $\tau$. Then the tautological commutative diagram

$$\begin{array}{ccc} B_1[\ell^r] & \to & \text{Hom}(A_1[\ell^r], \mu_{\ell^r}) \\ \downarrow \tau & & \downarrow \tilde{\sigma} \\ B_2[\ell^r] & \to & \text{Hom}(A_2[\ell^r], \mu_{\ell^r}) \end{array}$$

of $\mathbb{Z}_\ell[[G_{\mathbb{Q}_\ell}]]$-modules induces the commutative diagram of abelian groups

$$\begin{array}{ccc} H^1(\mathbb{Q}_\ell, B_1[\ell^r]) & \cong & \text{Hom}(H^1(\mathbb{Q}_\ell, A_1[\ell^r]), \mathbb{Z}/\ell^r\mathbb{Z}) \\ \downarrow \tau_{\ell,r} & & \downarrow \tilde{\sigma}_{\ell,r} \\ H^1(\mathbb{Q}_\ell, B_2[\ell^r]) & \cong & \text{Hom}(H^1(\mathbb{Q}_\ell, A_2[\ell^r]), \mathbb{Z}/\ell^r\mathbb{Z}) \end{array}$$

where the isomorphisms in the rows are from the Tate local duality for finite Galois representations. Applying the proof in the previous case to $B_1$ and $B_2$, we obtain

$$\tau_{\ell,r}(H^1_f(\mathbb{Q}_\ell, B_1, \ell^r)) = H^1_f(\mathbb{Q}_\ell, B_2, \ell^r). \tag{7}$$

By [**Wi**, Proposition 1.4.], $H^1_f(\mathbb{Q}_\ell, A_i, \ell^r)$ is the exact annihilator of $H^1_f(\mathbb{Q}_\ell, B_i, \ell^r)$ under the perfect pairing associated with the isomorphisms in the rows of the above commutative diagram, i.e., under the Tate duality. Combining this with the commutativity of the above diagram and Eq (7), we obtain $\sigma_{\ell,r}(H^1_f(\mathbb{Q}_\ell, A_1, \ell^r)) = H^1_f(\mathbb{Q}_\ell, A_2, \ell^r)$. Then the theorem again follows from Proposition 1.2. ∎

## 2.2 Kummer type congruences

We now compare the congruences in Theorem 1 with the classical Kummer congruences which states that for $\ell - 1 \nmid k$ and $k \equiv k' \mod (\ell - 1)\ell^{n-1}$,

$$(1 - \ell^{k-1})\frac{B_k}{k} \equiv (1 - \ell^{k'-1})\frac{B_{k'}}{k'} \mod \ell^n.$$

Define

$$H^1_f(\mathbb{Q}, (\mathbb{Q}/\mathbb{Z})(k)) = \bigoplus_\ell H^1_f(\mathbb{Q}, (\mathbb{Q}_\ell/\mathbb{Z}_\ell)(k)).$$

Bloch and Kato [**B-K**], applying the Iwasawa Main Conjecture verified by Mazur and Wiles [**M-W**], have proved that, for $k$ positive even,

$$\mid \frac{2(k-1)!}{(2\pi i)^k} \zeta(k) \mid = \frac{\mid H^1_f(\mathbb{Q}, (\mathbb{Q}/\mathbb{Z})(k)) \mid}{\mid H^0(\mathbb{Q}, (\mathbb{Q}/\mathbb{Z})(k)) \mid}$$

up to a power of 2. On the other hand, we have the well-known formula of Euler

$$\mid \frac{2(k-1)!}{(2\pi i)^k} \zeta(k) \mid = \mid \frac{B_k}{k} \mid.$$

We thus have

$$\frac{|H^1_f(\mathbb{Q},(\mathbb{Q}/\mathbb{Z})(k))|}{|H^0(\mathbb{Q},(\mathbb{Q}/\mathbb{Z})(k))|} = |\frac{B_k}{k}| \tag{8}$$

up to a power of 2. Using this equation, we prove

**Proposition 2.1** 1. Let $v_p$ be the p-valuation on $\mathbb{Q}$. Then

$$|H^0(\mathbb{Q},(\mathbb{Q}/\mathbb{Z})(k))| = \text{denominator of } |\frac{B_k}{k}| = \prod_{(p-1)|k} p^{v_p(k)+1}.$$

2. $|H^1_f(\mathbb{Q},(\mathbb{Q}/\mathbb{Z})(k))|$ equals to the numerator of $|\frac{B_k}{k}|$ up to a power of 2.

3. The greatest common divisor $(|H^1_f(\mathbb{Q},(\mathbb{Q}/\mathbb{Z})(k))|,|H^0(\mathbb{Q},(\mathbb{Q}/\mathbb{Z})(k))|)$ is a power of 2

**Proof:** The second statement follows directly from equation (8) and the first statement. Since the numerator and denominator of $\frac{B_k}{k}$ are relatively prime, the third statement follows from the first two statement. Therefore we only need to prove the first statement.

We first prove $|H^0(G_\mathbb{Q},(\mathbb{Q}/\mathbb{Z})(k))| = \prod_{(p-1)|k} p^{v_p(k)+1}$. From the isomorphism

$$(\mathbb{Q}/\mathbb{Z})(k) \cong \bigoplus_p (\mathbb{Q}_p/\mathbb{Z}_p)(k),$$

we only need to show that, for any positive integer $n$, $(\mathbb{Q}_p/\mathbb{Z}_p)(k)^{G_\mathbb{Q}} \geq p^n$ if and only if $(p-1)p^{n-1} \mid k$. Fix such an $n$. By definition, the action of $G_\mathbb{Q}$ on $(\mathbb{Q}_p/\mathbb{Z}_p)(1)[p^n] = \mu_{p^n}$ is given by the character $\chi_n = \chi \pmod{p^n} : G_\mathbb{Q} \to (\mathbb{Z}/p^n\mathbb{Z})^\times$ of order $(p-1)p^{n-1}$. Then the action of $G_\mathbb{Q}$ on $(\mathbb{Q}_p/\mathbb{Z}_p)(k)[p^n]$ is given by the character $\chi_n^k$. Therefore,

$$\begin{aligned}
|((\mathbb{Q}_p/\mathbb{Z}_p)(k))^{G_\mathbb{Q}}| &\geq p^n \\
\iff ((\mathbb{Q}_p/\mathbb{Z}_p)(k)[p^n])^{G_\mathbb{Q}} &= (\mathbb{Q}_p/\mathbb{Z}_p)(k)[p^n] \\
\iff \chi_n^k &= \text{id}\,|_{(\mathbb{Q}_p/\mathbb{Z}_p)(k)[p^n]} \\
\iff (p-1)p^{n-1} &\mid k.
\end{aligned}$$

To prove the equation

$$\text{denominator of } |\frac{B_k}{k}| = \prod_{(p-1)|k} p^{v_p(k)+1}$$

we apply the following two basic facts on Bernoulli numbers [**I-R**, **Wa**].

1. (Adam) If $p-1 \nmid k$, then $\frac{B_k}{k}$ is a $p$-integer.

2. (Carlitz) if $p$ is an odd prime and $p - 1 \mid k$, then $\frac{B_k + p^{-1} - 1}{k}$ is a $p$-integer.

Thus if $p - 1 \nmid k$, then $p$ does not divide the denominator of $B_k/k$. So the equation holds. If $p - 1 \mid k$, then by the formula of Carlitz, $\frac{B_k}{k} + \frac{p-1}{pk} = c$ for some $c \in \mathbb{Z}_p$. Thus by the trigonometry inequality in $\mathbb{Z}_p$,

$$v_p(\frac{B_k}{k}) = v_p(c - \frac{p-1}{pk}) = \min\{v_p(c), v_p(\frac{p-1}{pk})\} = v_p(\frac{p-1}{pk}).$$

So the $p$-power dividing the denominator of $B_k/k$ is the same as the $p$-power dividing the denominator of $(p^{-1} - 1)/k = (p-1)/pk$. This power is one higher than the $p$-power dividing $k$, so it is $v_p(k) + 1$. This proves the equation in the case when $p - 1 \mid k$. ∎

As an application of Theorem 1, we obtain the follow congruences which can be regarded as the $\ell$-part of the Kummer congruences.

**Theorem 2** *If $\ell - 1 \nmid k$ and $k \equiv k' \mod (\ell - 1)\ell^{n-1}$, then*

$$\ell - \text{part } (1 - \ell^{k-1})\frac{B_k}{k} \equiv \ell - \text{part } (1 - \ell^{k'-1})\frac{B_{k'}}{k'} \mod \ell^n.$$

**Proof:** From the definition,

$$H^1_f(\mathbb{Q}, (\mathbb{Q}_\ell/\mathbb{Z}_\ell)(k)) = H^1_f(\mathbb{Q}, (\mathbb{Q}/\mathbb{Z})(k))[\ell^\infty].$$

If $\ell - 1 \nmid k$, then from Proposition 2.1 we have $H^0(\mathbb{Q}, (\mathbb{Q}_\ell/\mathbb{Z}_\ell)(k)) = 0$. So

$$\begin{aligned}
\mid H^1_f(\mathbb{Q}, (\mathbb{Q}_\ell/\mathbb{Z}_\ell)(k)) \mid &= \ell - \text{part } \mid H^1_f(\mathbb{Q}, (\mathbb{Q}/\mathbb{Z})(k)) \mid \\
&= \ell - \text{part} \frac{\mid H^1_f(\mathbb{Q}, (\mathbb{Q}/\mathbb{Z})(k)) \mid}{\mid H^0(\mathbb{Q}, (\mathbb{Q}/\mathbb{Z})(k)) \mid} \\
&= \ell - \text{part } \mid \frac{B_k}{k} \mid \\
&= \ell - \text{part}(1 - \ell^{k-1}) \mid \frac{B_k}{k} \mid
\end{aligned}$$

Thus the congruences in the theorem follows from Theorem 1. ∎

## 3 Hecke characters

We now consider an analog of Theorem 1 for two dimensional representations. Let $K$ be an imaginary quadratic field with ring of integers $\mathcal{O}_K$. Let $E$ be an elliptic curve over $\mathbb{Q}$ with complex multiplication by $\mathcal{O}_K$. Let $\psi$ be the Hecke character associated to $E$ by the theory of complex multiplication. Let $\bar{\psi}$ be the complex conjugation of $\psi$.

For integers $k$ and $j$, define $\varphi = \varphi(k,j) = \psi^k \bar\psi^j$. For a finite prime $\lambda$ of $\mathcal{O}_K$, let $K_\lambda$ be the completion with the ring of integers $\mathcal{O}_\lambda$. Let

$$\varphi_\lambda : G_K \to \mathcal{O}_\lambda^\times$$

be the Weil realization of $\varphi$ at $\lambda$. Let $K_\lambda(\varphi_\lambda)$ (resp. $\mathcal{O}_\lambda(\varphi_\lambda)$) be the representation of $G_K$ with coefficient in $K_\lambda$ (resp. $\mathcal{O}_\lambda$). Now we use $K_\lambda(\varphi_\lambda)$ to construct a representation of $G_\mathbb{Q}$. Let $\ell$ be the prime of $\mathbb{Q}$ under $\lambda$, if $\ell$ splits in $K$ into $\lambda\lambda^*$, then $K_\lambda = \mathbb{Q}_\ell$. If $\ell$ is non-split in $K$, then $K_\lambda$ is of dimension two over $\mathbb{Q}_\ell$. In this case $\psi_\lambda$ and $\varphi_\lambda = \psi_\lambda^{k-j} \chi_\ell^j$ are representations of $G_\mathbb{Q}$. Define

$$M_\ell = M(k,j)_\ell = \begin{cases} \mathbb{Q}_\ell(\varphi_\lambda) \oplus \mathbb{Q}_\ell(\varphi_{\lambda^*}) & \text{if } p \text{ splits as } \lambda\lambda^* \text{ in } K \\ K_\lambda(\varphi_\lambda) & \text{if } p \text{ non-split, } \lambda | p \text{ in } K \end{cases}$$

Then $M_\ell$ is a representation of $G_\mathbb{Q}$ of dimension two over $\mathbb{Q}_\ell$. It is in fact the $\ell$-adic realization of the motive $M(k,j)$ for the algebraic Hecke character $\psi^k \bar\psi^j$ [**Gu3**]. Define a $G_\mathbb{Q}$-invariant $\mathbb{Z}_\ell$-lattice of $M(k,j)_\ell$ by

$$\mathcal{M}_\ell = \mathcal{M}(k,j)_\ell = \begin{cases} \mathbb{Z}_\ell(\varphi_\lambda) \oplus \mathbb{Z}_\ell(\varphi_{\lambda^*}) & \text{if } p \text{ splits as } \lambda\lambda^* \text{ in } K \\ \mathcal{O}_\lambda(\varphi_\lambda) & \text{if } p \text{ non-split, } \lambda | p \text{ in } K \end{cases}$$

Denote

$$A(k,j) = M(k,j)_\ell / \mathcal{M}(k,j)_\ell$$

with the induced representation of $G_\mathbb{Q}$. $A(k,j)$ is a cofree $\mathbb{Z}_\ell$-module of corank 2.

In this section we study the congruences between the Selmer groups of $A(k,j)$ for different pairs $(k,j)$. We will consider the cases when $\ell$ is split in $\mathcal{O}_K$ and non-split in $\mathcal{O}_K$ in the next two subsections. Some applications of these congruences to special values will also be given.

## 3.1 Case 1: $\ell$ splits

**Theorem 3** *Let $\ell$ be an odd prime of $\mathbb{Q}$ that splits in $K$. Let $(k_i, j_i)$, $i = 1, 2$ be two pairs of integers with $k_i > -j_i > 0$ and $H^0(\mathbb{Q}, A(k_1, j_1)) = 0$. If $A(k_1, j_1)[\ell^n] \cong A(k_2, j_2)[\ell^n]$ as $G_\mathbb{Q}$-representations, then $H^1_f(\mathbb{Q}, A(k_1, j_1))$ and $H^1_f(\mathbb{Q}, A(k_2, j_2))$ are algebraically congruent modulo $\ell^n$.*

**Proof:** Let $A = A(k,j)$ be either $A(k_1, j_1)$ or $A(k_2, j_2)$. It is known [**Gu3**] that the Selmer group $H^1_f(\mathbb{Q}, A)$ is the same as the Selmer groups of Greenberg [**Gr**], defined by

$$S_A^{\text{str}}(\mathbb{Q}) = \ker\left\{ H^1(\mathbb{Q}, A) \to \bigoplus_{p \neq \ell} \frac{H^1(\mathbb{Q}_p, A)}{H^1(g_p, A)_{\text{div}}} \oplus \frac{H^1(\mathbb{Q}_\ell, A)}{\text{im}(H^1(\mathbb{Q}_\ell, F^+ A))_{\text{div}}} \right\}.$$

Here

$$F^+ A = \mathbb{Z}_\ell(\varphi(k,j)_\lambda) \otimes (\mathbb{Q}_\ell / \mathbb{Z}_\ell).$$

Denote $N = N(k,j) = \mathbb{Z}_\ell(\varphi(k,j)_\lambda) \otimes (\mathbb{Q}_\ell/\mathbb{Z}_\ell)$. Then we further have
$$S^{\mathrm{str}}_{A(k,j)}(\mathbb{Q}) \cong S^{\mathrm{str}}_N(K),$$
where
$$S^{\mathrm{str}}_N(K) = \ker\left\{H^1(K,N) \to \bigoplus_{\mathfrak{p}\nmid\ell} \frac{H^1(K_\mathfrak{p},N)}{H^1(g_\mathfrak{p},N)_{\mathrm{div}}} \oplus \frac{H^1(K_\lambda,N)}{H^1(K_\lambda,N)_{\mathrm{div}}} \oplus H^1(K_{\lambda^*},N)\right\}.$$

We next recall more notations from [**Gu2**]. The reader is referred to there for more details and further references.

Let $\mathcal{K}_\infty = K(N)$. The Galois character $\varphi_\lambda$ is in fact given by the action of $\mathrm{Gal}(\bar{K}/K)$ on $N$
$$\varphi_\lambda = \psi^k \bar\psi_\lambda^{-j} : \mathrm{Gal}(\mathcal{K}_\infty/K) \to \mathrm{Aut}(N) \cong \mathbb{Z}_\ell^\times.$$

$\psi^k \bar\psi_\lambda^{-j}$ is a character of infinite order. Hence the image of $\varphi_\lambda$ is of finite index in $\mathbb{Z}_\ell^\times$. Thus $\mathrm{Gal}(\mathcal{K}_\infty/K)$ decomposes as $\mathrm{Gal}(\mathcal{K}_\infty/K) \cong \Delta \times \Gamma$ with $\Delta$ cyclic of order dividing $\ell - 1$ and $\Gamma \cong \mathbb{Z}_\ell$ as abelian groups. Also $\varphi$ decomposes accordingly into $\chi$ on $\Delta$ and $\kappa$ on $\Gamma$. Thus the fixed field $K_\infty$ of $\Delta$ is a $\mathbb{Z}_\ell$-extension of $K$. We will identify $\mathrm{Gal}(K_\infty/K)$ with $\Gamma$ by the natural isomorphism between them. Let $\mathcal{K}$ be the fixed field of $\Gamma$. Then $\mathcal{K} = K(N[\ell])$.

Define $F = K(E[\ell])$, $\Delta' = \mathrm{Gal}(F/\mathcal{K})$, $\delta = \mathrm{Gal}(F/K)$, $\mathcal{F}_\infty = K(E[\ell^\infty])$, $L_\infty = F(N)$. Let $F_\infty$ be the unique $\mathbb{Z}_\ell^2$-extension of $K$ contained in $\mathcal{F}_\infty$. Let $\Gamma_\ell = \mathrm{Gal}(F_\infty/K)$, $\Gamma' = \mathrm{Gal}(\mathcal{F}_\infty/L_\infty)$. Then we have $\mathrm{Gal}(\mathcal{F}_\infty/F) \cong \Gamma_\ell$ and $\mathrm{Gal}(\mathcal{F}_\infty/K) \cong \Gamma_\ell \times \delta$. In summary, we have the following diagram.

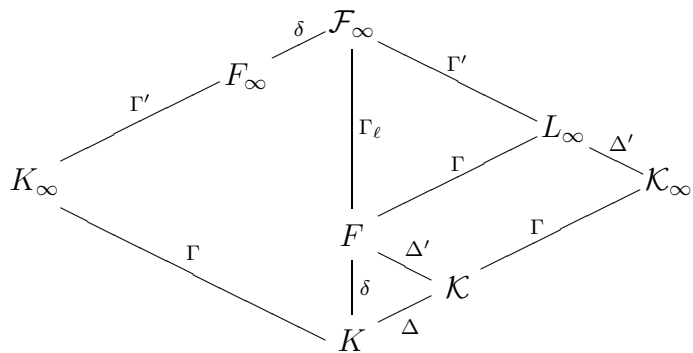

For any abelian extension $L$ of $K$, let $M(L)$ be the maximal abelian pro-$\ell$-extension of $L$ which is unramified outside the primes of $L$ over $\lambda$. Let $X(L)$ be $\mathrm{Gal}(M(L)/L)$. Combining Proposition 2.2 and Proposition 2.3 from [**Gu2**], we have
$$S^{\mathrm{str}}_N(K) \cong \mathrm{Hom}(X(L_\infty)^\chi, N)^\Gamma. \tag{9}$$

In the proof of Theorem 1 from [**Gu2**] we show that the restriction map $X(\mathcal{F}_\infty) \to X(L_\infty)$ induces an exact sequence of $\mathrm{Gal}(L_\infty/K)$-modules (note that $\Gamma'$ is denoted by $H$ there)
$$0 \to (X(\mathcal{F}_\infty))_{\Gamma'} \xrightarrow{\pi} X(L_\infty) \to \mathbb{Z}_\ell \to 0.$$

Since $\chi$ is a non-trivial character on $\delta$ by assumption, taking the $\chi$ component in the exact sequence, we get
$$(X(\mathcal{F}_\infty))_{\Gamma'}^\chi \cong X(L_\infty)^\chi.$$

Combining with Equation (9), we get
$$\begin{aligned} S_N^{\mathrm{str}}(K) &\cong \mathrm{Hom}((X(\mathcal{F}_\infty))_{\Gamma'}^\chi, N)^\Gamma \\ &\cong \mathrm{Hom}((X(\mathcal{F}_\infty))_{\Gamma'}, N)^{\mathrm{Gal}(L_\infty/K)} \\ &\cong (\mathrm{Hom}(X(\mathcal{F}_\infty), N)^{\mathrm{Gal}(L_\infty/K)})^{\Gamma'} \\ &\cong \mathrm{Hom}(X(\mathcal{F}_\infty), N)^{\mathrm{Gal}(\mathcal{F}_\infty/K)}. \end{aligned}$$

Now the theorem follows from Proposition 1.2. ∎

The congruences between Selmer groups obtained in Theorem 3 enables us to obtain congruences between special values of Hecke $L$-functions. This is accomplished by making use of the relation between Selmer groups and special values established by the Bloch-Kato conjecture for Hecke characters which has been verified in many cases [**Gu3**, **Han**]. Such congruences are similar to the congruences related to two variable $p$-adic $L$-functions [**Ya**].

**Theorem 4** *Let $\ell$ be an odd prime of $\mathbb{Q}$ that splits in $K$. Let $(k_i, j_i)$, $i = 1, 2$ be two pairs of integers with $k_i > -j_i > 0$ and $H^0(\mathbb{Q}, A(k_1, j_1)) = 0$. If $(k_1, j_1) \equiv (k_2, j_2) \pmod{(\ell - 1)\ell^{n-1}}$, then*
$$\ell-\mathrm{part}\left(\frac{2\pi}{\sqrt{d_K}}\right)^{j_1} \frac{L(\bar{\psi}^{j_1+k_1}, k_1)}{\Omega_\infty^{j_1+k_1}} \equiv \ell-\mathrm{part}\left(\frac{2\pi}{\sqrt{d_K}}\right)^{j_2} \frac{L(\bar{\psi}^{j_2+k_2}, k_2)}{\Omega_\infty^{j_2+k_2}} \pmod{\ell^n},$$

*where $d_K$ is the discriminant of $K$ and $\Omega_\infty$ is the period of the elliptic curve associated to $\psi$.*

**Proof:** By [**Gu2**, Theorem 1], for $i = 1, 2$,
$$\ell - \mathrm{part}\left(\frac{2\pi}{\sqrt{d_K}}\right)^{j_i} \Omega_\infty^{-(j_i+k_i)} L(\bar{\psi}^{j_i+k_i}, k_i) = |H_f^1(\mathbb{Q}, A(k_i, j_i))|.$$

Then the theorem follows from Theorem 3. ∎

## 3.2 Case 2: $\ell$ is non-split

We now consider the case when $\ell$ is inert in $K$. Then $\ell$ extends to a prime $\lambda$ in $K$. Recall that
$$A = A(k,j) = \mathcal{O}_\lambda(\varphi_\lambda) \otimes (K_\lambda/\mathcal{O}_\lambda).$$

Let $\mathcal{G} = \mathrm{Gal}(\mathcal{F}_\infty/K)$. Assume that $G_\lambda \stackrel{\mathrm{def}}{=} \mathrm{Gal}(\bar{K}_\lambda/K_\lambda)$ acts non-trivially on $A[\lambda]$. Then as in [**Ru**, Lemma 11.6], the inflation-restriction sequence and local classfield theory give an isomorphism

$$\phi : H^1(K_\lambda, A) \cong \mathrm{Hom}_\mathcal{G}(U(\mathcal{F}_\infty), A), \tag{10}$$

where
$$U(\mathcal{F}_\infty) = \varprojlim {}_n U(K(E[p^n]))$$

and $U(K(E[p^n]))$ is the group of local units of $K(E[p^n]) \otimes_K K_\lambda$ which are congruent to 1 modulo the primes above $\lambda$. From this isomorphism we get an injective $\mathcal{O}_\lambda$-module homomorphism

$$\eta = \eta_A : H^1_f(K_\lambda, A) \to \mathrm{Hom}_\mathcal{G}(U(\mathcal{F}_\infty), A). \tag{11}$$

Let $U_A$ be the annihilator of $\eta(H^1_f(K_\lambda, A))$ in $U(\mathcal{F}_\infty)$. For a positive integer $n$, we also have an injective $\mathcal{O}_\lambda/\lambda^n$-module homomorphism

$$\eta_n = \eta_{A,n} : H^1_f(K_\lambda, A)[\lambda^n] \to \mathrm{Hom}_\mathcal{G}(U(\mathcal{F}_\infty), A)[\lambda^n] \cong \mathrm{Hom}_\mathcal{G}(U(\mathcal{F}_\infty), A[\lambda^n]). \tag{12}$$

**Definition 3.1**  1. Denote $U_{A,n}$ for the annihilator of $\eta_n(H^1_f(K_\lambda, A)[\lambda^n])$ in $U(\mathcal{F}_\infty)$.

2. For a fixed pair $(k_0, j_0)$ with $k_0 > -j_0 \geq 0$ and a fixed positive integer $n$, denote $S(k_0, j_0, n)$ for the set of $A(k,j)$ with $k > -j \geq 0$ such that $A(k,j)[\lambda^n] \cong A(k_0, j_0)[\lambda^n]$ as $G_\mathbb{Q}$-modules.

3. For $A_1 = A(k_1, j_1)$, $A_2 = A(k_2, j_2) \in S(k_0, j_0, n)$, define $A_1$ and $A_2$ to be r-congruent modulo $\lambda^n$ if $U_{A_1,n} = U_{A_2,n}$.

We expect that any two representations from $S(k_0, j_0, n)$ are r-congruent to each other modulo $\lambda^n$. However its verification seems to require delicate properties of the Bloch-Kato exponential map for the $\lambda$-adic Hecke character $\varphi_\lambda$ associated to $A(k,j)$. We could nevertheless obtain the following weaker result.

**Proposition 3.1** *Assume that $G_\lambda$ acts non-trivially on $A(k_0, j_0)[\lambda]$. For almost all $A(k,j)$ in $S(k_0, j_0, n)$, there are infinitely many elements $A(k', j')$ in $S(k_0, j_0, n)$ such that $A(k,j)$ is r-congruent to $A(k', j')$ modulo $\lambda^n$.*

We first give the following result in algebra.

**Lemma 3.1** *Let $R$ and $S$ be two commutative rings. Let $A$ be a $R$-module, $B$ be an $S$-module and $C$ be a $R \otimes S$-module. There is a canonical isomorphism*

$$\xi : \mathrm{Hom}_R(A, \mathrm{Hom}_S(B, C)) \to \mathrm{Hom}_S(B, \mathrm{Hom}_R(A, C)).$$

**Proof:** For $f \in \mathrm{Hom}_R(A, \mathrm{Hom}_S(B, C))$, define $\xi(f)$ by

$$(\xi(f)(b))(a) = (f(a))(b), \quad \forall a \in A, \ b \in B.$$

$\xi(f)$ is clearly well-defined and is additive in both $a$ and $b$. Further, $\xi(f)(b)$ is $R$-linear since

$$(\xi(f)(b))(ra) = (f(ra))(b) = (rf(a))(b) = r((f(a))(b)) = r((\xi(f)(b))(a)).$$

$\xi(f)$ is $S$-linear since

$$\begin{aligned}(\xi(f)(sb))(a) &= f(a)(sb) = s(f(a)(b)) = s((\xi(f)(b))(a)) \\ &= (s(\xi(f)(b)))(a) = ((s\xi(f))(b))(a).\end{aligned}$$

By changing $A$ and $B$, we similarly obtain a homomorphism $\xi'$. It is easily verified that $\xi$ and $\xi'$ are the inverse of each other, thus proving that $\xi$ is an isomorphism. ∎

**Proof of Proposition 3.1:** Let $B_{\mathrm{dr},\lambda}$ be the ring of Fontaine with the filtration $\mathrm{Fil}^i B_{\mathrm{dr},\lambda}, i \in \mathbb{Z}$. It is known that

$$\mathrm{D}R(K_\lambda(\varphi_\lambda)) = H^0(K_\lambda, K_\lambda(\varphi) \otimes B_{\mathrm{dr},\lambda}) \cong K_\lambda.$$

From the assumption that $k > -j > 0$ we get

$$\mathrm{D}R^0(K_\lambda(\varphi_\lambda)) = H^0(K_\lambda, K_\lambda(\varphi) \otimes \mathrm{Fil}^0 B_{\mathrm{dr},\lambda}) = 0.$$

Since $G_{K_\lambda}$ acts non-trivially on $A[\lambda]$ and $A(k_0, j_0)[\lambda] \cong A(k, j)[\lambda]$, $G_\lambda$ acts non-trivially on $A(k, j)[\lambda]$. Then

$$H^0(K_\lambda, K_\lambda(\varphi_\lambda)) = 0.$$

So from [**B-K**, Corollary 3.8.4] we obtain $H^1_f(K_\lambda, K_\lambda(\varphi_\lambda)) \cong K_\lambda$. It follows that $H^1_f(K_\lambda, A)$ is a cofree $\mathcal{O}_\lambda$-module of corank one. Therefore $H^1_f(K_\lambda, A)[\lambda^n]$ is a free $\mathcal{O}_\lambda/\lambda^n$-module of rank one. By fixing an $\mathcal{O}_\lambda/\lambda^n$-module isomorphism

$$f : H^1_f(K_\lambda, A)[\lambda^n] \cong \mathcal{O}_\lambda/\lambda^n$$

and applying Lemma 3.1, we obtain the following sequence of isomorphisms

$$\begin{aligned}&\mathrm{Hom}_{\mathcal{O}_\lambda}(H^1_f(K_\lambda, A)[\lambda^n], \mathrm{Hom}_{\mathcal{G}}(U(\mathcal{F}_\infty), A)[\lambda^n]) \\ &\cong \mathrm{Hom}_{\mathcal{G}}(U(\mathcal{F}_\infty), \mathrm{Hom}_{\mathcal{O}_\lambda}(H^1_f(K_\lambda, A)[\lambda^n], A[\lambda^n])) \\ &\cong \mathrm{Hom}_{\mathcal{G}}(U(\mathcal{F}_\infty), \mathrm{Hom}_{\mathcal{O}_\lambda}(\mathcal{O}_\lambda/\lambda^n, A[\lambda^n])) \\ &\cong \mathrm{Hom}_{\mathcal{G}}(U(\mathcal{F}_\infty), A[\lambda^n])\end{aligned}$$

This way, $\eta_n$ corresponds to a $\delta_n = \delta_{A,n} \in \text{Hom}_{\mathcal{G}}(U(\mathcal{F}_\infty), A[\lambda^n])$. For another $\mathcal{O}_\lambda/\lambda^n$-module isomorphism $f' : H^1_f(K_\lambda, A)[\lambda^n] \cong \mathcal{O}_\lambda/\lambda^n$, we have $f' = \alpha \circ f$ where $\alpha$ is a automorphism of $\mathcal{O}_\lambda/\lambda^n$. Thus $\ker(\delta_n)$ does not depend on the choice of $f$ and is a well-defined $\Lambda$-submodule of $U(\mathcal{F}_\infty)$. From the definition of the isomorphism $\xi$ constructed in Lemma 3.1, it is easy to see that in fact $\ker(\delta_n) = U_n$.

Now let $A$ be an element in $S(k_0, j_0, n)$. So there is a $G_\mathbb{Q}$-module isomorphism

$$g : A[\lambda^n] \cong A_0[\lambda^n].$$

This isomorphism induces a $\mathcal{O}_\lambda/\lambda^n$-module isomorphism

$$\text{Hom}_{\mathcal{G}}(U(\mathcal{F}_\infty), A[\lambda^n]) \cong \text{Hom}_{\mathcal{G}}(U(\mathcal{F}_\infty), A_0[\lambda^n]),$$

under which $\delta_n$ is sent to $\tilde{\delta}_n = g \circ \delta_n$. Note that $\ker(\tilde{\delta}_n) = \ker(\delta_n) = U_n$ is independent of the choice of the isomorphism $g$ since any other $G_\mathbb{Q}$-module isomorphism

$$g' : A[\lambda^n] \cong A_0[\lambda^n]$$

is of the form $g' = \beta \circ g$ where $\beta$ is a $G_\mathbb{Q}$-module automorphism of $A_0[\lambda^n]$. Therefore, $A_1$ is r-congruent to $A_2$ modulo $\lambda^n$ if and only if $\ker(\tilde{\delta}_{A_1,n}) = \ker(\tilde{\delta}_{A_2,n})$. It is well-known [**Ru**] that $U(\mathcal{F}_\infty)$ is a free $\mathbb{Z}_\ell[[\mathcal{G}]]$-module of rank 2. Combining this with the fact that $A_0[\lambda^n]$ is finite, we see that the set $\text{Hom}_{\mathcal{G}}(U(\mathcal{F}_\infty), A_0[\lambda^n])$ is finite. Thus the set

$$T = \{\ker(\tilde{\delta}) \mid \tilde{\delta} \in \text{Hom}_{\mathcal{G}}(U(\mathcal{F}_\infty), A_0[\lambda^n])\}$$

is finite. Consequently, for almost all $A \in S(k_0, j_0, n)$, there are infinitely many $A' \in S(k_0, j_0, n)$ such that $\ker \tilde{A}_n = \ker \tilde{A}'_n$. This proves the proposition. ■

**Theorem 5** *Let $(k_i, j_i)$, $i = 1, 2$ be two pairs of integers with $k_i > -j_i > 0$. Let $\ell$ be an odd prime that is inert in $K$ and let $\lambda$ be the prime of $K$ over $\ell$. For a positive integer $n$, if $A(k_1, j_1)[\ell^n]$ and $A(k_2, j_2)[\ell^n]$ are isomorphic as $G_\mathbb{Q}$-representations, if $H^0(K_\lambda, A(k_1, j_1)) = 0$ and if $A(k_1, j_1)$ and $A(k_2, j_2)$ are r-congruent modulo $\lambda^n$, then $H^1_f(\mathbb{Q}, A(k_1, j_1))$ and $H^1_f(\mathbb{Q}, A(k_2, j_2))$ are algebraically congruent modulo $\lambda^n$.*

**Proof:** The assumption that $A(k_1, j_1)(K_\lambda) = 0$ is equivalent to the assumption that $G_\lambda = \text{Gal}(\bar{K}_\lambda/K_\lambda)$ acts non-trivially on $A(k_1, j_1)[\lambda]$. Since $A(k_1, j_1)[\lambda^n]$ and $A(k_2, j_2)[\lambda^n]$ are isomorphic $G_\mathbb{Q}$-representations, $G_\lambda$ acts non-trivially on $A(k_2, j_2)[\lambda]$. Now let $A = A(k, j)$ be either $A(k_1, j_1)$ or $A(k_2, j_2)$. Define $X(\mathcal{F}_\infty)$ as in the case when $\ell$ splits, considered in section 3.1. Let

$$H^1_f(K, A)' = \ker\left\{H^1(K, A) \to \bigoplus_{\mathfrak{p} \neq \ell} \frac{H^1(K_\mathfrak{p}, A)}{H^1_f(K_\mathfrak{p}, A)}\right\}.$$

Then we have the exact sequence

$$0 \to H^1_f(K, A) \to H^1_f(K, A)' \to \frac{H^1(K_\lambda, A)}{H^1_f(K_\lambda, A)}. \qquad (13)$$

Since $G_\lambda$ and therefore $G_K$ act non-trivially on $A[\lambda]$, the inflation-restriction sequence gives

$$H^1_f(K, A)' \cong \mathrm{Hom}_{\mathcal{G}}(X(\mathcal{F}_\infty), A). \qquad (14)$$

Further, we have the commutative diagram

$$\begin{array}{ccc} H^1(K, A) & \xrightarrow{\phi} & \mathrm{Hom}_{\mathcal{G}}(X(\mathcal{F}_\infty), A) \\ \downarrow & & \downarrow \\ H^1(K_\lambda, A) & \cong & \mathrm{Hom}_{\mathcal{G}}(U(\mathcal{G}_\infty), A) \end{array}$$

where the $\phi$ is the isomorphism from (10), the left vertical map is the cohomology restriction map and the right vertical map is the Artin map from classfield theory. Under the isomorphism (10), $H^1_f(K_\lambda, A)$ is identified with a $\mathcal{O}_\lambda$-submodule $H_A$ of $\mathrm{Hom}_{\mathcal{G}}(U(\mathcal{F}_\infty), A)$. Let $U_A$ be the annihilator of $H_A$ in $U(\mathcal{F}_\infty)$. Then $U_A$ is a $\mathcal{G}$-submodule of $U(\mathcal{F}_\infty)$ and we have $H_A \cong \mathrm{Hom}_{\mathcal{G}}(U(\mathcal{F}_\infty)/U_A, A)$. Thus from the exact sequence

$$0 \to \mathrm{Hom}_{\mathcal{G}}(U(\mathcal{F}_\infty)/U_A, A) \to \mathrm{Hom}_{\mathcal{G}}(U(\mathcal{F}_\infty), A) \to \mathrm{Hom}_{\mathcal{G}}(U_A, A)$$

we obtain the injection

$$\frac{H^1(K_\mathfrak{p}, A)}{H^1_f(K_\mathfrak{p}, A)} \hookrightarrow \mathrm{Hom}_{\mathcal{G}}(U_A, A).$$

From this injection and the isomorphism (14), the sequence (13) can be identified with the exact sequence

$$0 \to H^1_f(K, A) \to \mathrm{Hom}_{\mathcal{G}}(X(\mathcal{F}_\infty), A) \to \mathrm{Hom}_{\mathcal{G}}(U_A, A), \qquad (15)$$

where the second map comes from the Artin map. From this exact sequence we have the following commutative diagram with exact rows

$$\begin{array}{ccccccc} 0 & \to & H^1_f(K, A) & \to & \mathrm{Hom}_{\mathcal{G}}(X(\mathcal{F}_\infty), A) & \to & \mathrm{Hom}_{\mathcal{G}}(U_A, A) \\ & & \downarrow \lambda^n & & \downarrow \lambda^n & & \downarrow \lambda^n \\ 0 & \to & H^1_f(K, A) & \to & \mathrm{Hom}_{\mathcal{G}}(X(\mathcal{F}_\infty), A) & \to & \mathrm{Hom}_{\mathcal{G}}(U_A, A) \end{array}$$

Applying the snake lemma gives

$$0 \to H^1_f(K, A)[\lambda^n] \to \mathrm{Hom}_{\mathcal{G}}(X(\mathcal{F}_\infty), A)[\lambda^n] \to \mathrm{Hom}_{\mathcal{G}}(U_A, A)[\lambda^n].$$

Equivalently,

$$0 \to H^1_f(K, A)[\lambda^n] \to \mathrm{Hom}_{\mathcal{G}}(X(\mathcal{F}_\infty), A[\lambda^n]) \to \mathrm{Hom}_{\mathcal{G}}(U_A, A[\lambda^n]). \qquad (16)$$

Now given $A_1 = A(k_1, j_1)$ and $A_2 = A(k_2, j_2)$ as in the theorem, Since they are r-congruent modulo $\lambda^n$, we have $U_{A_1} = U_{A_2}$. Then a $G_\mathbb{Q}$-module isomorphism $A_1[\lambda^n] \cong A_2[\lambda^n]$ gives the following commutative diagram

$$\begin{array}{ccc} \mathrm{Hom}_\mathcal{G}(X(\mathcal{F}_\infty), A_1[\lambda^n]) & \to & \mathrm{Hom}_\mathcal{G}(U_{A_1}, A_1[\lambda^n]) \\ \downarrow & & \downarrow \\ \mathrm{Hom}_\mathcal{G}(X(\mathcal{F}_\infty), A_2[\lambda^n]) & \to & \mathrm{Hom}_\mathcal{G}(U_{A_2}, A_2[\lambda^n]) \end{array}$$

in which both of the column maps are isomorphisms. It follows from (16) that

$$H^1_f(K, A_1)[\lambda^n] \cong H^1_f(K, A_2)[\lambda^n].$$

Taking the $\mathrm{Gal}(K/\mathbb{Q})$-invariants, and noting that $\ell$ is odd by assumption, we have

$$H^1_f(\mathbb{Q}, A_1)[\lambda^n] \cong H^1_f(\mathbb{Q}, A_2)[\lambda^n].$$

∎

## 4 Adjoint representations

In this section, let $K$ be a number field and let $K_\lambda$ be the completion of $K$ at a finite prime $\lambda$ with the ring of integers $\mathcal{O}_\lambda$. Let $\ell$ be the prime of $\mathbb{Q}$ under $\lambda$. Let $M$ be a $G_\mathbb{Q}$-representation that is cofree over $\mathcal{O}_\lambda$ with corank two. Define $A = A(M) = \cup_n \mathrm{End}(M[\lambda^n])$, $A^0 = A^0(M) = \cup_n \mathrm{End}^0(M[\lambda^n])$ with the induced $G_\mathbb{Q}$-actions. Here

$$\mathrm{End}^0(M[\lambda^n]) = \{\text{trace zero endomorphisms of } N[\lambda^n]\}.$$

**Theorem 6** *Let $V_i$, $i = 1, 2$, be $\lambda$-adic representations of $G_\mathbb{Q}$ of dimension two that are unramified for almost all primes $p \neq \ell$ and crystalline at $\ell$ with Hodge-Tate type $(0, k-1)$, for some $2 \leq k < \ell$. Let $T_i$ be $G_\mathbb{Q}$-invariant lattices of $V_i$. Let $M_i = V_i/T_i$. Assume that $H^0(\mathbb{Q}, A^0(M_i))$ is finite and $H^0(I_p, A(M_i))$ is $\lambda$-divisible for $p \neq \ell$ and $i = 1, 2$. If $M_1[\lambda^n] \cong M_2[\lambda^n]$ as $G_\mathbb{Q}$-modules, then*

$$H^1_f(\mathbb{Q}, A(M_1))[\lambda^n] \cong H^1_f(\mathbb{Q}, A(M_2))[\lambda^n],$$

$$H^1_f(\mathbb{Q}, A^0(M_1))[\lambda^n] \cong H^1_f(\mathbb{Q}, A^0(M_2))[\lambda^n].$$

Before giving the proof of the theorem, we first display a corollary. Let $f$ be a newform of weight $k \geq 2$, level $N$ and with coefficients in $K$. For a fixed finite prime $\lambda$ of $K$, let $T(f) = T_\lambda(f)$ be the $\mathcal{O}_\lambda[[G_\mathbb{Q}]]$-module associated to $f$ by Deligne. Let $M(f) = M_\lambda(f) = T \otimes_{\mathcal{O}_\lambda} (K_\lambda/\mathcal{O}_\lambda)$. It is known that if $\lambda \nmid N$, then $V(f) = T(f) \otimes_{\mathcal{O}_\lambda} K_\lambda$ is unramified at almost all primes $p \neq \ell$ and is crystalline at $\lambda$. Assume that $V(f)$ is minimally ramified in the sense of [**Di**], then $H^0(I_p, A(M(f)))$ is $\lambda$-divisible for all $p \neq \ell$. If we further assume that $V(f)$ is absolutely irreducible, then $H^0(\mathbb{Q}, \mathrm{End}(V(f))) = \mathrm{End}_{G_\mathbb{Q}}(V)$ has only scalars. Since $\mathrm{End}(V(f)) = E_\lambda \oplus \mathrm{End}^0(V(f))$, we necessarily have $H^0(\mathbb{Q}, \mathrm{End}^0(V(f))) = 0$, and hence $H^0(\mathbb{Q}, A^0(M(f)))$ is finite. We have thus obtained

**Corollary 4.1** *Let $f_i$, $i = 1, 2$ be two newforms of weight $k$, level $N_i$. Let $\lambda$ be a finite prime of $K$ such that the prime $\ell$ of $\mathbb{Q}$ under $\lambda$ is greater than $k$, $\lambda$ does not divide $N_1 N_2$. Assume that $V(f_i)$ is absolutely irreducible, and is minimally ramified at $p \neq \ell$. If $M(f_1)[\lambda^n]$ is isomorphism to $M(f_2)[\lambda^n]$ as $G_\mathbb{Q}$-modules, then*

$$H^1_f(\mathbb{Q}, A(M(f_1)))[\lambda^n] \cong H^1_f(\mathbb{Q}, A(M(f_2)))[\lambda^n],$$

$$H^1_f(\mathbb{Q}, A^0(M(f_1)))[\lambda^n] \cong H^1_f(\mathbb{Q}, A^0(M(f_2)))[\lambda^n].$$

**Remark:** It can be shown that, for almost all $\lambda$, $V = V_\lambda$ is absolutely irreducible and is minimally ramified at $p \neq \ell$ [**DFG**].

**Proof of Theorem 6:** Denote $A_i$ for $A(M_i)$ and $A_i^0$ for $A(M_i^0)$, $i = 1, 2$. Let $A$ be either $A_1$ or $A_2$ and let $A^0$ be either $A_1^0$ or $A_2^0$. Note that $A(M) \cong A^0 \oplus K_\lambda/\mathcal{O}_\lambda$. Thus $H^1_f(\mathbb{Q}, A(M)) \cong H^1_f(\mathbb{Q}, A^0(M)) \oplus H^1_f(\mathbb{Q}, K_\lambda/\mathcal{O}_\lambda)$. Since [**Fl**]

$$H^1_f(\mathbb{Q}, \mathbb{Q}_\ell/\mathbb{Z}_\ell) \cong \text{Hom}(\text{Gal}(H/\mathbb{Q}), \mathbb{Q}_\ell/\mathbb{Z}_\ell)$$

where $H$ is the Hilbert class field of $\mathbb{Q}$, we have $H^1_f(\mathbb{Q}, \mathbb{Q}_\ell/\mathbb{Z}_\ell) = 0$ and consequently

$$H^1_f(\mathbb{Q}, K_\lambda/\mathcal{O}_\lambda) \cong H^1_f(\mathbb{Q}, \mathbb{Q}_\ell/\mathbb{Z}_\ell) \otimes_{\mathbb{Z}_\ell} \mathcal{O}_\lambda = 0.$$

Therefore we only need to prove the first equation in the theorem.

By the second statement of Proposition 1.2, we only need to prove that, for $1 \leq r \leq n$,

$$|H^1_f(F, A_1, \lambda^r)| = |H^1_f(F, A_2, \lambda^r)|, \tag{17}$$

For this purpose, we compare $H^1_f(\mathbb{Q}, A)$ with another version of the Selmer group of $A$, first defined by Wiles [**Wi**].

Again let $A$ be either $A_1$ or $A_2$. Let $\ell$ be the prime of $\mathbb{Q}$ below $\lambda$. Fix $r \geq 1$. Fontaine and Laffaille [**F-L**, **DDT**] have constructed a category $\mathcal{MF}_r^{f,k}$ of filtered $\mathcal{O}_\lambda$-modules of finite length, annihilated by $\lambda^r$, and with filtration length $k$, and shown that there is a fully faithful covariant functor from $\mathcal{MF}_r^{f,k}$ to a subcategory $\mathcal{FG}_r^k$ of the category $\mathfrak{Mod}_r$ of $(\mathcal{O}_\lambda/\lambda^r)[G_\ell]$-modules, thereby gives an anti-equivalence of the two categories.

Define $\text{Ext}^1_{\text{crys}}(M[\lambda^r], N[\lambda^r])$ to be the subgroup of $\text{Ext}^1(M[\lambda^r], N[\lambda^r])$ consisting of extensions from the subcategory $\mathcal{FG}_r^k$. Taking the inverse image of this subgroup under the canonical isomorphism

$$H^1(\mathbb{Q}_\ell, A[\lambda^r]) \cong \text{Ext}^1(M[\lambda^r], N[\lambda^r]),$$

we get a subgroup of $H^1(\mathbb{Q}_\ell, A[\lambda^r])$ which we denote by $H^1_w(\mathbb{Q}_\ell, A[\lambda^r])$. Such a subgroup is first studied by Wiles [**Wi**] in the case when $k = 2$.

For $r \geq 1$, the sequence of subgroups $\{H^1_w(\mathbb{Q}_\ell, A[\lambda^r])\}_r$ form a subsystem of the direct system $\{H^1(\mathbb{Q}_\ell, A[\lambda^r])\}_r$. So it makes sense to define

$$H^1_w(\mathbb{Q}_\ell, A) = \varinjlim H^1_w(\mathbb{Q}_\ell, A[\lambda^r]) \subseteq H^1(\mathbb{Q}_\ell, A)$$

When $\lambda \nmid p$, define

$$H^1_w(\mathbb{Q}_p, A[\lambda^r]) = H^1(G_p/I_p, (A[\lambda^r])^{I_p}),$$

and define

$$H^1_w(\mathbb{Q}_p, A) = H^1(G_p/I_p, A^{I_p}).$$

Then we have

$$H^1_w(\mathbb{Q}_p, A) = \varinjlim H^1_w(\mathbb{Q}_p, A[\lambda^r]).$$

Define Wiles' Selmer group of $A$ by

$$H^1_w(\mathbb{Q}, A) = \ker\left\{ H^1(\mathbb{Q}, A) \to \bigoplus_{p < \infty} \frac{H^1(\mathbb{Q}_p, A)}{H^1_w(\mathbb{Q}_p, A)} \right\}.$$

Similarly define

$$H^1_w(\mathbb{Q}, A[\lambda^r]) = \ker\left\{ H^1(\mathbb{Q}, A[\lambda^r]) \to \bigoplus_{p < \infty} \frac{H^1(\mathbb{Q}_p, A[\lambda^r])}{H^1_w(\mathbb{Q}_p, A[\lambda^r])} \right\}.$$

We again have

$$H^1_w(\mathbb{Q}, A) = \varinjlim H^1_w(\mathbb{Q}, A[\lambda^r]).$$

By definition, if $A_1[\lambda^r]$ and $A_2[\lambda]$ are isomorphic as $G_\mathbb{Q}$-modules, then we have the abelian group isomorphism

$$H^1_w(\mathbb{Q}, A_1[\lambda^r]) \cong H^1_w(\mathbb{Q}, A_2[\lambda^r]).$$

Therefore to prove Equation (17) and hence the theorem, we only need to prove

$$H^1_f(\mathbb{Q}, A, \lambda^r) = H^1_w(\mathbb{Q}, A[\lambda^n]).$$

For this, we only need to prove

**Lemma 4.1** *For any finite prime $p$ of $\mathbb{Q}$, we have*

$$H^1_f(\mathbb{Q}_p, A, \lambda^r) = H^1_w(\mathbb{Q}_p, A[\lambda^r]).$$

**Proof:** First consider the case when $p \neq \ell$. By assumption, $A^{I_p} = H^0(I_p, A)$ is $\lambda$-divisible. Thus we have the exact sequence

$$0 \to A^{I_p}[\lambda^r] \to A^{I_p} \xrightarrow{\lambda^r} A^{I_p} \to 0$$

of $g_p$-modules, where $g_p = \text{Gal}(\bar{\mathbb{Q}}_p/\mathbb{Q}_p)/I_p$. Taking the Galois cohomology sequence from this sequence and comparing, through the inflation map, with the Galois cohomology sequence from the short exact sequence of $G_p$-modules

$$0 \to A[\lambda^r] \xrightarrow{j_r} A \xrightarrow{\lambda^r} A \to 0,$$

we get the following commutative diagram with exact rows

$$\begin{array}{ccccccccc}
0 & \to & H^0(\mathbb{Q}_p, A)/\lambda^r H^0(\mathbb{Q}_p, A) & \to & H^1(g_p, (A^{I_p})[\lambda^r]) & \to & H^1(g_p, A^{I_p})[\lambda^r] & \to & 0 \\
& & \| & & \downarrow & & \downarrow & & \\
0 & \to & H^0(\mathbb{Q}_p, A)/\lambda^r H^0(\mathbb{Q}_p, A) & \to & H^1(\mathbb{Q}_p, A[\lambda^r]) & \xrightarrow{j_r} & H^1(\mathbb{Q}_p, A)[\lambda^r] & \to & 0
\end{array}$$

where the vertical arrows are injective. By definition

$$H_w^1(\mathbb{Q}_p, A[\lambda^r]) = H^1(g_p, (A[\lambda^r])^{I_p}) = H^1(g_p, (A^{I_p})[\lambda^r]).$$

It also follows from the definition that $H_f^1(\mathbb{Q}_p, A)$ is the maximal divisible subgroup of $H^1(g_p, A^{I_p})$. Since $A^{I_p}$ is divisible, it follows that $H^1(g_p, A^{I_p}) \cong A^{I_p}/(\gamma - 1)A^{I_p}$ is also divisible. Here $\gamma$ is a topological generator of $g_p \cong \hat{\mathbb{Z}}$. Therefore

$$H_f^1(\mathbb{Q}_p, A) = H^1(g_p, A^{I_p}).$$

Then from the above commutative diagram we get

$$H_f^1(\mathbb{Q}_p, A, \lambda^r) \stackrel{\text{def}}{=} j_r^{-1}(H_f^1(\mathbb{Q}_p, A)) = H^1(g_p, (A^{I_p})[\lambda^r]) = H_w^1(\mathbb{Q}_p, A[\lambda]).$$

The proof in the case when $p = \ell$ follows the idea of [**Wi**, Proposition 1.3]. But note that $H_f^1$ there is $H_w^1$ here and $H_F^1$ there is $H_f^1$ here. For more details, see [**DFG**]. Adapting Wiles' argument and using the explicit descriptions of $H_w^1(\mathbb{Q}_\ell, A[\lambda^n])$ [**F-L**, **Ra**, **DDT**], we obtain

$$H_w^1(\mathbb{Q}_p, A) = H_f^1(\mathbb{Q}_p, A).$$

In particular they are $\lambda$-divisible with the same $\mathcal{O}_\lambda$-corank, say $c$. By the definition of $H_f^1(\mathbb{Q}_p, A, \lambda^n)$, we have

$$H_w^1(\mathbb{Q}_p, A[\lambda^n]) \subseteq H_f^1(\mathbb{Q}_p, A, \lambda^n)$$

and

$$|H_f^1(\mathbb{Q}_p, A, \lambda^n)| = |(\mathcal{O}_\lambda/\lambda^n)^c| \cdot |H^0(\mathbb{Q}_p, A[\lambda^n])|.$$

However this last number is also the order of $H_w^1(\mathbb{Q}_p, A[\lambda^n])$ [**DDT**, 2.5]. Therefore the two groups must be the same. This proves the Lemma when $p = \ell$. ∎

Now the proof of Theorem 6 is completed. ∎

# References

[**B-K**] S. Bloch and K. Kato, *L-functions and Tamagawa numbers of motives,* The Grothendieck Festschrift, Vol 1. Birkhäuser(1990), 333–400.

[**DDT**] H. Darmon, F. Diamond, R. Taylor, *Fermat's Last Theorem,* In: Current Development in Mathematics, 1995, International Press, 1-154.

[**Di**] F. Diamond, *An extension of Wiles' results,* In: Modular forms and Fermat's last theorem, Springer-Verlag, 1997.

[**DFG**] F. Diamond, M. Flach, L. Guo, *On the Bloch-Kato conjecture for adjoint motives of modular forms,* in preparation.

[**Fl**] M. Flach, *A generalization of the Cassels-Tate pairing,* J. Reine Angew. Math. **412** (1990), 113-127.

[**F-L**] J.-M. Fontaine and G. Laffaille, *Construction de représentations p-adiques,* Ann. Sci. Ec. Norm. Super. **15** (1982), 547-608.

[**Gr**] R. Greenberg, *Iwasawa theory for p-adic representations,* Advanced Studies in Pure Mathematics **17**, Academic Press (1989), 97–137.

[**Gu1**] L. Guo, *On a generalization of Tate dualities with application to Iwasawa theory,* Comp. Math. **85** (1993), 125-161.

[**Gu2**] L. Guo, *General Selmer groups and critical value of Hecke L-functions,* Math. Ann. **297** (1993), 221-233.

[**Gu3**] L. Guo, *On the Bloch-Kato conjecture for Hecke L-functions,* J. Number Theory **57** (1996), 340-365.

[**Han**] B. Han, *On Bloch-Kato conjecture of Tamagawa numbers for Hecke characters of imaginary quadratic number fields,* Ph.D. thesis, University of Chicago, 1997.

[**I-R**] K. Ireland and M. Rosen, *A classical introduction to modern number theory,* 2nd ed. (1990), Springer-Verlag.

[**Iw**] K. Iwasawa, *On $\mathbb{Z}_\ell$-extensions of algebraic number fields,* Ann. of Math.. **98** (1973), 146-326.

[**Ma**] B. Mazur, *Deforming Galois representations,* In: Galois groups over $\mathbb{Q}$, Y. Ihara, K Ribet, J-P. Serre, ed., MSRI Publ. **16**, Springer-Verlag, New York, Berlin, Heidelberg, 1989, pp. 385-437.

[**M-W**] B.Mazur and A. Wiles, *Class fields of abelian extensions of Q,* Invent. Math. **76** (1984), 179-330.

[**Ra**] R. Ramakrishna, *On a variation of Mazur's deformation functor,* Compositio Math. **87**(1993), 269-286.

[**Ru**] K. Rubin, *The "main conjecture" of Iwasawa theory for imaginary quadratic fields,* Invent. Math. **103** (1991), 25–68.

[**Wa**] L. Washington, *Introduction to cyclotomic fields,* 2nd edition, (1996), Springer-Verlag.

[**Wi**] A. Wiles, *Modular elliptic curves and Fermat's Last Theorem,* Annals of Math. **141** (1995), 443-551.

[**Ya**] R.Yager, *On two variable p-adic L-functions,* Ann. Math. **115** (1982), 411-449.